\documentclass[a4paper,12pt,centertags,psamsfonts]{amsart}
\usepackage{amsmath,latexsym}
\usepackage{amssymb}
\usepackage{times}
\DeclareSymbolFont{SY}{U}{psy}{m}{n}
\DeclareMathSymbol{\emptyset}{\mathord}{SY}{'306}
\chardef\bslash=`\\ 

\newtheorem*{thm}{Theorem} 

\newtheorem{prop}{Proposition}

\theoremstyle{definition}
\newtheorem{defn}{Definition}[]
\theoremstyle{remark}
\newtheorem*{rem}{Remark} 
\newtheorem*{notation}{Notation}

\newcommand{\propref}[1]{Proposition~\ref{#1}}

\newcommand{\Na}{\mathbb{N}}
\newcommand{\fC}{\mathfrak{C}}

\newcommand{\Zi}{\mathbb{Z}}

\newcommand{\R}{\mathbb{R}}
\newcommand{\fR}{\mathfrak{R}}

\newcommand{\C}{\mathbb{C}}

\let\bsy\boldsymbol
\def\P{\bsy P}
\def\T{\bsy T}
\def\t{\bsy t}

\def\p{\bsy p}
\def\f{\bsy f}
\def\K{\bsy K}
\def\k{\bsy k}

\def\F{\bsy F}

\def\H{\mathcal{H}}

\newcommand{\Ss}{\mathcal{S}}

\newcommand{\wt}{\widetilde}

\DeclareMathOperator{\sgn}{sign}

\newcommand{\eval}[2][\right]{\relax
  \ifx#1\right\relax \left.\fi#2#1\rvert}

\let\abs=\envert

\let\norm=\enVert

\let\hnorm=\henVert
\newcommand{\kskob}[1]{\left[#1\right]}
\let\ksk=\kskob

\let\sk=\skob
\newcommand{\skoba}[1]{\left\langle#1\right\rangle}
\let\ska=\skoba

\let\hska=\hskoba

\let\krsk=\krskoba

\let\sets=\fskoba
\begin{document}
\renewcommand{\sectionmark}[1]{}
\title[Integral representations of
closed operators]{Integral representations of
closed operators \\ as bi-Carleman operators \\ with arbitrarily smooth kernels}
\author[I. M. Novitski\u i]{Igor M. Novitski\u i}
\address{Institute for Applied Mathematics, Russian Academy of Sciences,
92, Zaparina   Street, Khabarovsk 680 000, Russia}
\email{novim@iam.khv.ru}
\thanks{Research supported in part by grant N 03-1-0-01-009 from
the Far-Eastern Branch of the Russian Academy of Sciences. This paper
was written in November 2003, when the author enjoyed the hospitality of the
Mathematical Institute of Friedrich-Schiller-University, Jena, Germany}
\begin{abstract}
In this paper, we characterize all closed linear operators
in a separable Hilbert space which are unitarily equivalent to an integral
bi-Carleman operator in $L_2(\R)$ with bounded and \textit{arbitrarily}
smooth kernel on $\R^2$. In addition, we give an explicit
construction of corresponding unitary operators. The main result
is a qualitative sharpening of an earlier result of \cite{nov:91}.
\end{abstract}
\maketitle

Throughout, $\H$  will denote a separable Hilbert space
with the inner product $\hska{\cdot,\cdot}$ and the norm $\hnorm{\cdot}$,
$\fC(\H)$ the set of all closed linear operators densely defined in $\H$,
and $\C$, and $\Na$, and $\Zi$, the complex plane, the set of all positive
integers, the set of all integers, respectively.
For an operator $S$ in $\fC(\H)$, $S^*$ will denote the Hilbert space
adjoint of $S$.

An operator $S:D_S\to\H$ of $\fC(\H)$ is said to belong to the set
$\fC_{00}(\H)$ if there exist a  linear manifold $D$ dense in $\H$
and an orthonormal sequence $\{e_n\}$ such that
\begin{equation}\label{biCarleman}
\{e_n\}\subset D\subset D_S\cap D_{S^*},\quad
   \lim\limits_{n\to\infty}\hnorm{Se_n}=0,\quad
   \lim\limits_{n\to\infty}\hnorm{S^* e_n}=0.
\end{equation}

Let $\R$ be the real line $(-\infty,+\infty)$ equipped with the Lebesgue
measure, and let $L_2=L_2(\R)$ be the Hilbert space of (equivalence classes
of) measurable complex--valued functions on $\R$ equipped with the inner
product
$$
\ska{f,g}=\int_{\R} f(s)\overline{g(s)}\,ds
$$
and the norm
$\Vert f\Vert=\ska{f,f}^{\frac{1}2}$.
A linear operator
                 $T:D_T\to L_2$,
                 where the domain
                       $D_T$
                       is a dense linear manifold in
                                                    $L_2$,
is said to be \textit{integral\/} if there exists a measurable function
$\T$ on $\R^2$, a \textit{kernel\/}, such that, for every $f\in D_T$,
$$
               (Tf)(s)=\int_\R\T(s,t)f(t)\,dt
$$
for almost every $s$ in $\R$. A kernel $\T$ on $\R^2$ is said to be
\textit{Carleman\/} if $\T(s,\cdot) \in L_2$ for almost every fixed $s$ in
$\R$. An integral operator with a kernel $\T$ is called \textit{Carleman\/}
if $\T$ is a Carleman kernel, and it is called \textit{bi-Carleman\/} if both
$\T$ and $\T^*$ ($\T^*(s,t)=\overline{\T(t,s)}$) are Carleman kernels.
Every Carleman kernel, $\T$, induces a \textit{Carleman
function\/} $\t$ from $\R$ to $L_2$ by $\t(s)=\overline{\T(s,\cdot)}$
for all $s$ in ${\R}$ for which $\T(s,\cdot)\in L_2$.

The characterization of closed linear operators whose unitary orbits
contain integral operators was first studied by Neumann in \cite{Neu}
and by now is well understood \cite{Kor:book1}. The main results related to
this problem are concerned with the unitary equivalence of linear operators
to Carleman or bi-Carleman operators with measurable kernels.
Here is one of them \cite[p.~145]{Kor:book1}.
\begin{prop}\label{Kor}
A necessary and sufficient condition that an operator $S\in\fC(\H)$ be
unitarily equivalent to a bi-Carleman operator is that $S$ belong to
$\fC_{00}(\H)$.
\end{prop}

However, the characterization problem can also be formulated in terms of
kernels that satisfy various additional conditions. For example, given any
non-negative integer $m$, we consider the following question:
which operators are unitarily equivalent to a bi-Carleman operator with a
kernel $\K$ satisfying the conditions:
\begin{enumerate}
\renewcommand{\labelenumi}{(\roman{enumi})}
\item the function $\K$ and all its partial derivatives
on $\R^2$ up to order $m$ are in $C(\R^2,\C)$,
\item the Carleman function $\k$,
$\k(s)=\overline{\K(s,\cdot)}$,
and all its (strong) derivatives on ${\R}$ up to order $m$
are in $C(\R,L_2)$,
\item the conjugate transpose function $\K^*$,
$\K^*(s,t)=\overline{\K(t,s)}$,
satisfies Condition (ii), that is,
the Carleman function $\k^*$, $\k^*(s)
=\overline{\K^*(s,\cdot)}$, and all its (strong) derivatives on $\R$ up
to order $m$ are in $C(\R,L_2)$?
\end{enumerate}
Here and throughout  $C(X,B)$, where $B$ is a Banach space (with norm
$\norm{\cdot}_B$), denote the Banach space (with the norm $\norm{f}_{C(X,B)}
=\sup\limits_{x\in X}\,\norm{f(x)}_B$) of continuous $B$-valued functions
defined on a locally compact space $X$ and \textit{vanishing at infinity\/}
(that is, given any $f \in C(X,B)$ and $\varepsilon>0$, there exists a
compact subset $X(\varepsilon,f) \subset X$ such that
$\norm{f(x)}_{B}<\varepsilon$ whenever $x\not\in X(\varepsilon,f)$).

A function $\K$ that satisfies Conditions (i), (ii) is called a
\textit{$SK^m$-kernel\/} \cite{nov:Isra}. In addition, a $SK^m$-kernel $\K$
is called a \textit{$K^m$-kernel\/} (\cite{Nov:Lon}, \cite{nov:91}) if it
satisfies Condition (iii).

The next result is an answer to the above question; it gives a
characterization of closed linear operators representable as bi-Carleman
operators  with $K^m$-kernels (cf. \propref{Kor}).
\begin{prop}[\cite{nov:91},\cite{nov:92}]\label{msmooth}
Let $m$ be a fixed non-negative integer, and let $S\in\fC_{00}(\H)$.
Then there exists a unitary operator $U_m:\H\to L_2$ such that
$T=U_mSU_m^{-1}$ is a bi-Carleman operator having $K^m$-kernel.
\end{prop}
The purpose of this paper is to restrict the conclusion of \propref{msmooth}
to \textit{arbitrarily\/} smooth kernels. Now we define these kernels.

\begin{defn}
We say  that a function $\K$ is a \textit{$K^\infty$$(SK^\infty)$-kernel}
if it is a $K^m$($SK^m$)-kernel for each non-negative integer $m$
({\cite{Je}, \cite{SP}}).
\end{defn}

The following theorem is the main result of the present paper.

\begin{thm}\label{thm-main}
If $S\in\fC_{00}(\H)$,
then there exists a unitary operator $U_\infty:\H\to L_2$ such that the
operator $T=U_\infty SU^{-1}_\infty$ is a bi-Carleman operator having
$K^\infty$-kernel.
\end{thm}

\section{Proof of Theorem}
The proof has an algorithmic nature and consists of three steps.
In the first step we split the operator $S\in\fC_{00}(\H)$ in
order to construct auxiliary operators $J$, $\wt J$, $B$, $\wt B$,
$Q$, $\wt Q$. Using these operators, in the second step we
describe suitable orthonormal bases $\{u_n\}$ in $L_2$ and $\{f_n\}$ in
$\H$, and then we use the mentioned bases to construct a unitary
operator from $\H$ to $L_2$ which sends the basis $\sets{f_n}$
onto the basis $\sets{u_n}$, and fulfills the role of the operator
$U_\infty$ cited in the theorem. The remaining part of the proof
(step 3) is a direct verification that the constructed unitary
operator has the desired properties. Thus, the proof of the
theorem not only establishes the unitary equivalence itself, but
also indicates an explicit construction of $U_\infty$.

\subsection*{Step 1.} Let $S\in\fC_{00}(\H)$.
Assume, with no loss of generality,
that the sequence
$\sets{e_k}_{k=1}^\infty\subset D$ in
\eqref{biCarleman} satisfies the condition
\begin{equation}\label{null-seq}
\sum_k\krsk{\hnorm{Se_k}^{\frac14}+\hnorm{S^*e_k}^{\frac14}}\le1
\end{equation}
(the sum notation $\sum\limits_k$ will always be used instead of
the more detailed symbol $\sum\limits_{k=1}^\infty$).
Let $H$ denote a subspace spanned by the $e_k$'s, and let
$H^\perp$ be the orthogonal complement of $H$ in $\H$.
Since $S$, $S^*\in\fC(\H)$, we have
\begin{equation}\label{hvd}
H\subset D_S,\quad H\subset D_{S^*}.
\end{equation}
If $E$ is the orthogonal projection onto $H$, consider the decompositions
\begin{equation}\label{splitting}
  S=(1-E)S+ES, \quad S^*=(1-E)S^*+ES^*.
\end{equation}
Since $E\in\fR(\H)$, it follows via definition of the adjoint that
$(ES)^*=S^*E$, and $(ES^*)^*=SE$.
Observe that the operators
\begin{equation}\label{E}
J=SE,\quad \wt J=S^*E
\end{equation}
are nuclear because of \eqref{null-seq} and \eqref{hvd}.
 Assume, with no loss of generality, that
$\dim H^\perp=\infty$, and choose an orthonormal basis
$\sets{e_k^\perp}_{k=1}^\infty$ for $H^\perp$ so that
\begin{equation}\label{hkvd}
  \sets{e_k^\perp}_{k=1}^\infty\subset (1-E)D\subset D_{S}\cap D_{S^*}.
\end{equation}
For each $f\in D_{S}\cap D_{S^*}$ and for each $h\in \H$, let
\begin{equation}\label{defdknums}
\begin{gathered}
z(f)=\hnorm{Sf}+\hnorm{S^*f}, \\
d(h)=\hnorm{Jh}^\frac{1}4+\hnorm{J^*h}^\frac{1}4
     +\hnorm{\wt Jh}^\frac{1}4+\hnorm{\wt J^*h}^\frac{1}4.
\end{gathered}
\end{equation}
If
$$
J=\sum\limits_n s_n\hska{\cdot,p_n}q_n\ \text{and}\
\wt J=\sum\limits_n\wt s_n\hska{\cdot,\wt p_n}\wt q_n
$$
are the Schmidt decompositions
for $J$ and $\wt J$, respectively,
then the closedness of both $S$ and $S^*$ implies that,
for all $f\in D_{S^*}$ and $g\in D_S$,
\begin{equation}\label{ES}
\begin{gathered}
ES^*f=(ES^*)^{**}f=(SE)^*f=J^*f=\sum\limits_n  s_n\hska{g, q_n} p_n,\\
ESg=(ES)^{**}g=(S^*E)^*g=\wt J^*g=\sum\limits_n \wt s_n\hska{g,\wt q_n}\wt p_n;
\end{gathered}
\end{equation}
here the $s_n$ are the singular values of $J$ (eigenvalues of
$\krsk{JJ^*}^{\frac12}$), $\sets{p_n}$, $\sets{q_n}$  are
orthonormal sets (the $p_n$ are eigenvectors for $J^*J$ and the $q_n$ are
eigenvectors for $JJ^*$).

Define
\begin{equation}\label{B}
B=\sum_n s_{n}^{\frac1{4}}\hska{\cdot,p_n} q_{n},\quad
\wt B=\sum_n\wt s_{n}^{\frac1{4}}\hska{\cdot,\wt p_n}\wt q_{n},
\end{equation}
and observe that, by the Schwarz inequality,
\begin{equation}\label{schw}
\begin{gathered}
\hnorm{Bf}=\hnorm{\left(J^*J\right)^{\frac18}f}\leq\hnorm{J f}^{\frac14},\\
\hnorm{B^*f}=\hnorm{\left(JJ^*\right)^{\frac18}f}\leq \hnorm{J^*f}^{\frac14},
\\
\hnorm{\wt Bf}=\hnorm{\left(\wt J^*\wt J\right)^{\frac18}f}\leq\hnorm{\wt J f}^{\frac14},\\
\hnorm{\wt B^*f}=\hnorm{\left(\wt J\wt J^*\right)^{\frac18}f}\leq \hnorm{\wt J^*f}^{\frac14}
\end{gathered}
\end{equation}
if $\norm{f}=1$.
The operators $B$ and $\wt B$ play only an auxiliary role in what follows.

Define operators $Q=(1-E)S$, $\wt Q=(1-E)S^*$.
The property \eqref{hkvd}
guarantees that the following representations hold:
\begin{equation}\label{qoper}
\begin{gathered}
Qf=\sum\limits_k\hska{Qf,e_k^\perp} e_k^\perp=
  \sum\limits_k\hska{f,S^*e_k^\perp} e_k^\perp\quad\text{for all $f\in D_S$,}\\
\wt Qg=\sum\limits_k\hska{\wt Qg,e_k^\perp} e_k^\perp=
  \sum\limits_k\hska{g,Se_k^\perp} e_k^\perp\quad\text{for all $g\in D_{S^*}$}.
\end{gathered}
\end{equation}

\subsection*{Step 2.} This step is to  construct a candidate for the desired
unitary operator $U_\infty$ in the theorem.

\begin{notation}
If an equivalence class
$f\in L_2$ contains a function belonging to $C(\R,\C)$, then we shall use
$\ksk{f}$ to denote that function.
\end{notation}
Take orthonormal bases $\sets{f_n}$ for $\H$ and $\sets{u_n}$ for $L_2$
which satisfy the conditions:
\begin{enumerate}
\renewcommand{\labelenumi}{(\alph{enumi})}
\item the terms of the sequence
$\left\{\ksk{u_n} ^{(i)}\right\}$ of derivatives are in $C(\R,\C)$, for each $i$
(here and throughout, the letter $i$ is reserved for all non-negative integers),
\item for each $i$,
\begin{gather}
\begin{gathered}\label{hki}
\sum_k H_{k,i}<\infty,\ \sum_k z\sk{v_k}H_{m(k),i}<\infty,\
\sum_k z\sk{e_k^\perp}H_{n(k),i}<\infty
\\
\text{with $H_{k,i}=\norm{\ksk{h_k}^{(i)}}_{C(\R,\C)}$}\quad (k\in\Na),
\end{gathered}\\
\label{zndn}
\sum_k d(x_k)\left(G_{k,i}+1\right)<\infty\quad
\text{with $G_{k,i}=\norm{\ksk{g_k}^{(i)}}_{C(\R,\C)}$}\quad (k\in\Na),
\end{gather}
where $\sets{n(k)}_{k=1}^\infty$ and $\sets{m(k)}_{k=1}^\infty$ are
subsequences of $\Na$ such that
$\sets{m(k)}_{k=1}^\infty=\Na\setminus\sets{n(k)}_{k=1}^\infty$, and
$\sets{g_k}$, $\sets{h_k}$,
$\sets{x_k}$, and
$\sets{v_k}$, are orthonormal sets such that
\begin{equation}\label{SET}
\begin{gathered}
\{u_n\}=\{g_k\}_{k=1}^\infty\cup\{h_k\}_{k=1}^\infty,\quad
\{g_k\}_{k=1}^\infty\cap\{h_k\}_{k=1}^\infty=\varnothing,\\
\sets{f_n}=\sets{x_k}_{k=1}^\infty
\cup\sets{v_k}_{k=1}^\infty
\cup\sets{e_k^\perp}_{k=1}^\infty,
\quad
\sets{x_k}_{k=1}^\infty\subset\sets{e_k}_{k=1}^\infty,\\
\sets{v_k}_{k=1}^\infty=\sets{e_k}_{k=1}^\infty\setminus\sets{x_k}_{k=1}^\infty.
\end{gathered}
\end{equation}
\end{enumerate}

\begin{rem}
Let $\{ u_n\}$ be an orthonormal basis for $L_2$ such that, for
each $i$,
\begin{gather}\label{1}
\kskob{u_n}^{(i)}\in C(\R,\C)\quad(n\in\Na),\\
\label{2}\norm{\kskob{u_n}^{(i)}}_{C(\R,\C)}\le D_nA_i\quad(n\in\Na),\\
\label{3}
\sum_kD_{n_k}<\infty,
\end{gather}
where $\{D_n\}_{n=1}^\infty$, $\{A_i\}_{i=0}^\infty$ are sequences
of positive numbers, and $\{n_k\}_{k=1}^\infty$ is a subsequence
of $\Na$ such that
$\Na\setminus\sets{n_k}_{k=1}^\infty$ is a countable set.
Since $d(e_k)\to0$ and $z(v_k)\to0$ as $k\to\infty$, the basis $\{ u_n\}$
satisfies Conditions (a) and (b), with $h_k=u_{n_k}$
($k\in\Na$) and
$\{g_k\}_{k=1}^\infty=\{u_n\}\setminus\{h_k\}_{k=1}^\infty$.

To construct an example of
such basis $\{u_n\}$, consider a Lemari\'e-Meyer wavelet,
$$
u(s)=\dfrac1{2\pi}\int_{\R}e^{i\xi(\frac12+s)}
\sgn\xi b(|\xi|)\,d\xi\quad (s\in\R),
$$
with the bell function $b$ belonging to
$C^\infty(\R)$ (for construction of the Lemari\'e-Meyer wavelets
we refer to \cite{LeMe}, \cite[\S~4]{Ausch}, \cite[Example D, p.~62]{Her}).
Then $u$ belongs to the Schwartz class $\Ss(\R)$, and hence
all the derivatives $\kskob{u}^{(i)}$
are in $C(\R,\C)$.
The ``mother function'' $u$ generates an orthonormal basis for $L_2$ by
$$
u_{jk}(s)=2^{\frac j2}u(2^js-k)\quad  (j,\,k\in\Zi).
$$
 Rearrange, in a completely arbitrary manner, the orthonormal set
$\{u_{jk}\}_{j,\,k\in\Zi}$ into a simple sequence,
so that it becomes $\{ u_n\}_{n\in\Na}$. Since, in view of this rearrangement,
to each $n\in\Na$ there corresponds a unique pair of integers
$j_n$, $k_n$, and
conversely, we can write, for each $i$,
$$
\norm{\kskob{u_n}^{(i)}}_{C(\R,\C)}=\norm{\kskob{u_{j_nk_n}}^{(i)}}_{C(\R,\C)}
\le D_nA_i,
$$
where
$$
D_n=\begin{cases}
2^{j_n^2}&\text{if $j_n>0$,}\\
\left(\dfrac1{\sqrt{2}}\right)^{\abs{j_n}}&\text{if $j_n\le0$,}
\end{cases}
\qquad A_i=2^{\left(i+\frac12\right)^2}\norm{\kskob{u}^{(i)}}_{C(\R,\C)}.
$$
Whence it follows that if $\{n_k\}_{k=1}^\infty\subset\Na$ is a subsequence
such that $j_{n_k}\to -\infty$ as $k\to\infty$, then
$$\sum_kD_{n_k}<\infty.$$
Thus,  the basis $\{u_n\}$ satisfies Conditions \eqref{1} through \eqref{3} and,
consequently, Conditions (a) and (b).
\end{rem}

Let us return to the proof of the theorem.
Define a candidate for the desired unitary operator
$U_\infty:\H\to L_2$ on the basis vectors as follows:
\begin{equation}\label{uaction}
U_\infty x_k=g_k,\quad
U_\infty v_k=h_{m(k)},\quad U_\infty e_k^\perp=h_{n(k)}\quad \text{for all $k\in \Na$},
\end{equation}
where $\sets{n(k)}_{k=1}^\infty$ and $\sets{m(k)}_{k=1}^\infty$ are just
those sequences which occur in Condition \eqref{hki},
in the harmless assumption that, for each $k\in\Na$,
\begin{equation}\label{ektohnk}
 U_\infty f_k=u_k,\quad  y_k=U_\infty^{-1} h_k.
\end{equation}

\subsection*{Step 3.}
Show that the unitary operator $U_\infty$ defined in \eqref{uaction} has the
desired properties, that is, that $T= U_\infty S U_\infty ^{-1}$
is a bi-Carleman operator having $K^\infty$-kernel. For this purpose,
verify first that all the operators
$P= U_\infty Q U_\infty ^{-1}$,
$\wt P= U_\infty \wt Q U_\infty ^{-1}$,
$F= U_\infty J^* U_\infty ^{-1}$,
$\wt F= U_\infty \wt J^* U_\infty ^{-1}$ are Carleman operators having
$SK^\infty$-kernels.

 Using  \eqref{qoper}, \eqref{uaction}, one can write
\begin{equation} \label{P}
\begin{gathered}
Pf=\sum_k \ska{f,T^*h_{n(k)}} h_{n(k)}\quad\text{for all $f\in D_T= U_\infty D_S$},\\
\wt Pf=\sum_k \ska{g,Th_{n(k)}} h_{n(k)}\quad\text{for all $g\in D_{T^*}= U_\infty D_{S^*}$},
\end{gathered}
\end{equation}
where
\begin{equation}\label{Thk}
\begin{gathered}
T^*h_{n(k)}=\sum\limits_n\hska{e_k^\perp,Sf_n}u_n,\\
Th_{n(k)}=\sum\limits_n\hska{e_k^\perp,S^*f_n} u_n\quad (k\in\Na).
\end{gathered}
\end{equation}
Prove that, for any fixed $i$, the series
\begin{equation*}\label{imtuko} 
  \sum\limits_n\hska{ e_k^\perp, Sf_n} \ksk{u_n}^{(i)}(s), \quad
\sum\limits_n\hska{ e_k^\perp,S^*f_n} \ksk{u_n}^{(i)}(s) \quad
   (k\in\Na)
\end{equation*}
converge in the norm of $C(\R,\C)$.
Indeed, all these series are pointwise dominated on $\R$
by one series
$$
  \sum\limits_n\krsk{\hnorm{Sf_n}+\hnorm{S^*f_n}}\abs{\ksk{u_n}^{(i)}(s)},
$$
which converges uniformly on $\R$ because its component subseries
(see \eqref{uaction}, \eqref{SET})
\begin{gather*}
  \sum\limits_k\krsk{\hnorm{Jx_k}+\hnorm{\wt Jx_k}}\abs{\ksk{g_k}^{(i)}(s)},\\
  \sum\limits_k\krsk{\hnorm{Sv_k}+\hnorm{S^*v_k}}\abs{\ksk{h_{m(k)}}^{(i)}(s)},\\
  \sum\limits_k\krsk{\hnorm{Se_k^\perp}+\hnorm{S^*e_k^\perp}}\abs{\ksk{h_{n(k)}}^{(i)}(s)}
\end{gather*}
are in turn dominated by the convergent series
\begin{equation*}\label{domser} 
 \sum\limits_k d(x_k)G_{k,i}, \quad
 \sum\limits_k z(v_k)H_{m(k),i},\quad
 \sum\limits_k z(e_k^\perp)H_{n(k),i},
\end{equation*}
respectively (see \eqref{null-seq}, \eqref{E},  \eqref{defdknums}, \eqref{zndn}, \eqref{hki}).
Whence it follows  
that, for each $k\in\Na$,
\begin{equation}\label{supest}
\norm{\ksk{T^*h_{n(k)}}^{(i)}}_{C(\R,\C)}\le C^*_i,\quad
\norm{\ksk{Th_{n(k)}}^{(i)}}_{C(\R,\C)}\le C_i,
\end{equation}
with constants $C^*_i$ and $C_i$ independent of $k$.
From \eqref{defdknums} it follows also that
\begin{equation}\label{sest}
\norm{T^*h_{n(k)}}\le z\krsk{e_k^\perp},\quad
\norm{Th_{n(k)}}\le z\krsk{e_k^\perp}\quad(k\in\Na),
\end{equation}
since $U_\infty$ is unitary.
Consider functions $\P$, $\wt\P:\R^2\to\C$, and $\p$, $\wt\p:\R\to L_2$, defined, for all
$s$, $t\in\R$, by
\begin{equation} \label{Pp}
\begin{gathered}
\P(s,t)=\sum\limits_k\ksk{h_{n(k)}}(s)
         \overline{\ksk{T^*h_{n(k)}}(t)},\\
        \wt \P(s,t)=\sum\limits_k\ksk{h_{n(k)}}(s)
         \overline{\ksk{Th_{n(k)}}(t)},\\
\p(s)=\overline{\P(s,\cdot)}=\sum\limits_k
\overline{\ksk{h_{n(k)}}(s)}T^*h_{n(k)},\\
\wt\p(s)=\overline{\wt\P(s,\cdot)}=\sum\limits_k
\overline{\ksk{h_{n(k)}}(s)}Th_{n(k)}.
\end{gathered}
\end{equation}
The termwise differentiation theorem implies
that, for each $i$ and each non-negative integer $j$,
\begin{gather*}
\dfrac{\partial^{i+j}\P}{\partial s^i\partial t^j}(s,t)
=\sum\limits_{k}\ksk{h_{n(k)}}^{(i)}(s)
                             \overline{\ksk{T^*h_{n(k)}}^{(j)}(t)},\\
\dfrac{\partial^{i+j}\wt\P}{\partial s^i\partial t^j}(s,t)
=\sum\limits_{k}\ksk{h_{n(k)}}^{(i)}(s)
                             \overline{\ksk{Th_{n(k)}}^{(j)}(t)},
                             \\
\dfrac{d^i\p}{ds^i}(s)=
\sum\limits_k\overline{\ksk{h_{n(k)}}^{(i)}(s)}T^*h_{n(k)},\\
\dfrac{d^i\wt\p}{ds^i}(s)=
\sum\limits_k\overline{\ksk{h_{n(k)}}^{(i)}(s)}Th_{n(k)},
\end{gather*}
since, by \eqref{supest}, \eqref{sest}, and \eqref{hki},
the displayed series converge (absolutely) in $C(\R^2,\C)$, $C(\R,L_2)$, respectively.
Thus,
$$
\dfrac{\partial^{i+j}\P}{\partial s^i\partial t^j},\
\dfrac{\partial^{i+j}\wt\P}{\partial s^i\partial t^j}\in C(\R^2,\C),\quad
\text{and}\quad \dfrac{d^i\p}{ds^i},\
\dfrac{d^i\wt\p}{ds^i}\in C(\R,L_2).
$$
Observe also that, by \eqref{sest}, \eqref{hki}, and \eqref{Pp}, the series \eqref{P}
(viewed, of course, as ones with terms belonging to $C(\R,\C)$) converge
(absolutely) in $C(\R,\C)$-norm  to the functions
$$
\begin{gathered}
\kskob{Pf}(s)\equiv\skoba{f,\p(s)}\equiv\int_{\R}\P(s,t)f(t)\,dt,\\
\kskob{\wt Pg}(s)\equiv\skoba{g,\wt\p(s)}\equiv\int_{\R}\wt\P(s,t)g(t)\,dt,
\end{gathered}
$$
respectively. Thus, both $P:D_T\to L_2$ and $\wt P:D_{T^*}\to L_2$ are  Carleman operators, and $\P$ and $\wt\P$
are their $SK^\infty$-kernels,
respectively

Since, by \eqref{E}, $\hnorm{Se_k}=\hnorm{Je_k}$ and $\hnorm{S^*e_k}=\hnorm{\wt Je_k}$
for all $k$,
from \eqref{null-seq} it follows
via \eqref{schw} that the operators $B$ and $\wt B$ defined in \eqref{B} are
 nuclear, and hence
\begin{equation}\label{snums}
\sum_n s_{n}^{\frac1{2}}<\infty,\quad \sum_n\wt s_{n}^{\frac1{2}}<\infty.
\end{equation}
Then, according to \eqref{ES} and \eqref{B}, kernels which induce
the nuclear operators $F$ and $\wt F$ can be represented by the series
\begin{equation}\label{fkernels}
\sum_n s_n^{\frac12} U_\infty B^*q_n(s)
\overline{ U_\infty Bp_n(t)},\quad
\sum_n \wt s_n^{\frac12} U_\infty \wt B^*\wt q_n(s)
\overline{ U_\infty\wt  B\wt p_n(t)}
\end{equation}
convergent almost everywhere in $\R^2$. The functions used
in these bilinear expansions can be written as the series convergent in $L_2$:
$$
\begin{gathered}
   U_\infty Bp_k=\sum\limits_n\hska{ p_k,B^*f_n} u_n, \quad
   U_\infty B^*q_k=\sum\limits_n\hska{q_k,Bf_n} u_n,\\
   U_\infty \wt B\wt p_k=\sum\limits_n\hska{\wt p_k,\wt B^*f_n} u_n, \quad
   U_\infty \wt B^*\wt q_k=\sum\limits_n\hska{\wt q_k,\wt Bf_n} u_n\ (k\in\Na).
\end{gathered}
$$
Show that, for any fixed $i$, the functions
$\ksk{ U_\infty Bp_k}^{(i)}$, $\ksk{ U_\infty B^*q_k}^{(i)}$,
$\ksk{ U_\infty \wt B \wt p_k}^{(i)}$, $\ksk{ U_\infty\wt  B^*\wt q_k}^{(i)}$ ($k\in\Na$)
make sense,
are all in $C(\R,\C)$, and their $C(\R,\C)$-norms
are bounded independent of $k$. Indeed, all the series
\begin{equation*} 
\begin{gathered}
  \sum\limits_n\hska{p_k,B^*f_n}\ksk{u_n}^{(i)}(s),\quad
  \sum\limits_n\hska{q_k,Bf_n}\ksk{u_n}^{(i)}(s),\\
  \sum\limits_n\hska{\wt p_k,\wt B^*f_n}\ksk{u_n}^{(i)}(s),\quad
  \sum\limits_n\hska{\wt q_k,\wt Bf_n}\ksk{u_n}^{(i)}(s)\quad (k\in\Na)
\end{gathered}
\end{equation*}
are dominated by one series
$$
  \sum\limits_n\krsk{\hnorm{B^*f_n}+\hnorm{Bf_n}+
  \hnorm{\wt B^*f_n}+\hnorm{\wt Bf_n}}
  \abs{\ksk{u_n}^{(i)}(s)}.
$$
This series converges uniformly in $\R$, since it consists of two
uniformly convergent in $\R$ subseries (see \eqref{uaction}, \eqref{ektohnk})
\begin{equation*}
\begin{gathered}
  \sum\limits_k\krsk{\hnorm{B^*x_k}+\hnorm{Bx_k}
  +\hnorm{\wt B^*x_k}+\hnorm{\wt Bx_k}}\abs{\ksk{g_k}^{(i)}(s)},\\
  \sum\limits_k\krsk{\hnorm{B^*y_k}+\hnorm{By_k}
  +\hnorm{\wt B^*y_k}+\hnorm{\wt By_k}}
  \abs{\ksk{h_k}^{(i)}(s)},
\end{gathered}
\end{equation*}
which are dominated by the following convergent series
$$
\sum\limits_k d(x_k)G_{k,i}, \quad
\sum\limits_k 2\krsk{\norm{B}+\norm{\wt B}}H_{k,i},
$$
respectively
(see \eqref{E}, \eqref{null-seq}, \eqref{defdknums}, \eqref{schw}, \eqref{zndn}, \eqref{hki}).
Thus, for  functions $\F$, $\wt\F:\R^2\to \C$, $\f$, $\wt\f:\R\to L_2$,
defined by
\begin{equation*}
\begin{gathered}
\F(s,t)=\sum_n s_n^{\frac12}\ksk{ U_\infty B^*q_n}(s)\overline{\ksk{ U_\infty Bp_n}(t)},\\
\wt\F(s,t)=\sum_n \wt s_n^{\frac12}\ksk{ U_\infty \wt B^*\wt q_n}(s)\overline{\ksk{ U_\infty \wt B\wt p_n}(t)},\\
\f(s)=\overline{\F(s,\cdot)}=\sum_n s_n^{\frac12}\overline{\ksk{ U_\infty B^*q_n}(s)} U_\infty Bp_n,\\
\wt \f(s)=\overline{\wt \F(s,\cdot)}=\sum_n\wt  s_n^{\frac12}\overline{\ksk{ U_\infty\wt  B^*\wt q_n}(s)} U_\infty\wt B\wt p_n,
\end{gathered}
\end{equation*}
one can write, for all non-negative integers $i$, $j$ and all $s$, $t\in\R$,
\begin{equation*}
\begin{gathered}
\dfrac{\partial^{i+j}\F}{\partial s^i\partial t^j}(s,t)
=\sum_n s_n^{\frac12}\ksk{ U_\infty B^*q_n}^{(i)}(s)\overline{\ksk{ U_\infty Bp_n}^{(j)}(t)},\\
\dfrac{\partial^{i+j}\wt \F}{\partial s^i\partial t^j}(s,t)
=\sum_n\wt  s_n^{\frac12}\ksk{ U_\infty\wt  B^*\wt q_n}^{(i)}(s)\overline{\ksk{ U_\infty\wt  B\wt p_n}^{(j)}(t)},\\
\dfrac{d^i\f}{ds^i}(s)
=\sum_n s_n^{\frac12}\overline{\ksk{ U_\infty B^*q_n}^{(i)}(s)} U_\infty Bp_n,\\
\dfrac{d^i\wt\f}{ds^i}(s)
=\sum_n \wt s_n^{\frac12}\overline{\ksk{ U_\infty\wt  B^*\wt q_n}^{(i)}(s)} U_\infty \wt B\wt p_n,
\end{gathered}
\end{equation*}
where  the series converge in $C(\R^2,\C)$, $C(\R,L_2)$, respectively,
because of \eqref{snums}. This implies that $\F$ and $\wt\F$
are $SK^\infty$-kernels of $F$ and $\wt F$, respectively.

Using \eqref{splitting} and \eqref{ES}, we obtain decompositions
$$
  Tf=Pf+\wt Ff \quad (f\in D_T), \quad
  T^*f=\wt Pf+Ff \quad (f\in D_T^*)
$$
all of whose terms are Carleman operators. Therefore $T$ and $T^*$ are
also Carleman, and their kernels $\K$ and $\K^*$
can be defined by
$$
  \K(s,t)= \P(s,t)+\wt \F(s,t), \quad
  \K^*(s,t)=\wt\P(s,t)+\F(s,t)\quad(s,t \in\R),
$$
respectively, where all the terms are $SK^\infty$-kernels.
Moreover, we have $\K(s,t)=\overline{\K^*(t,s)}$, $\K(\cdot, t)=
\overline{\K^*(t, \cdot)}$ for all $s$, $t\in \R$
(cf. \cite[Corollary~IV.2.17]{Kor:book1}). Thus,
the operator $T=U_\infty SU_\infty^{-1}$ is a bi-Carleman operator,
with the kernel $\K$, which is a $K^\infty$-kernel.
The proof of the theorem is complete.

\section*{Acknowledgments} The author thanks the Mathematical Institute
of the  University of Jena for its hospitality, and specially
W.~Sickel and H.-J.~Schmei\ss er for useful remarks and fruitful
discussion on applying wavelets in integral representation theory.

\bibliographystyle{amsplain}

\enddocument